\documentclass[10pt,journal,onecolumn]{IEEEtran}

\usepackage{bm, amssymb, amsmath, graphics, mathtools}
\newcommand{\T}{^{\mbox{\tiny T}}}
\newcommand{\B}[1]{{\bm #1}}
\newcommand{\ds}{\displaystyle}
\newcommand{\dd}{\; \text{d}}

\begin{document}

\title{Least-squares Solutions  \protect\\ of Linear Differential Equations}

\author{Daniele Mortari \\ \vspace{0.5cm} \emph{dedicated to John Lee Junkins}\footnote{Professor, Aerospace Engineering, Texas A\&M University, College Station, TX $77843$-$3141$, USA. E-mail: mortari$@$tamu.edu}}

\maketitle

\begin{abstract}
    This study shows how to obtain \emph{least-squares solutions} to initial and boundary value problems to nonhomogeneous linear differential equations with nonconstant coefficients of \emph{any} order. However, without loss of generality, the approach has been applied to second order differential equations. The proposed method has two steps. The first step consists of writing a \emph{constrained expression}, introduced in Ref. \cite{Mortari}, that has embedded the differential equation constraints. These expressions are given in term of a new unknown function, $g (t)$, and they satisfy the constraints, no matter what $g (t)$ is. The second step consists of expressing $g (t)$ as a linear combination of $m$ independent known basis functions, $g (t) = \B{\xi}\T \B{h} (t)$. Specifically, Chebyshev orthogonal polynomials of the first kind are adopted for the basis functions. This choice requires rewriting the differential equation and the constraints in term of a new independent variable, $x\in[-1, +1]$. The procedure leads to a set of linear equations in terms of the unknown coefficients vector, $\B{\xi},$ that is then computed by \emph{least-squares}. Numerical examples are provided to quantify the solutions accuracy for initial and boundary values problems as well as for a control-type problem, where the state is defined in one point and the costate in another point.
\end{abstract}

\noindent {\bf Acronyms used throughout this paper}
\begin{tabbing}
  DE ~~ \= $\to$ \= Differential Equation \\
  IVP \> $\to$ \> Initial Value Problem \\
  BVP \> $\to$ \> Boundary Value Problem \\
  LS \> $\to$ \> Least-Squares
\end{tabbing}

\section{Introduction}

The $n$-th order nonhomogeneous ordinary linear Differential Equation (DE) with nonconstant coefficients is the equation
\begin{equation}\label{e01}
    \ds\sum_{i = 0}^n f_i (t) \, \dfrac{\dd^i y (t)}{\dd t^i} = f (t),
\end{equation}
where $f (t)$ and the $n + 1$ functions, $f_i (t)$, can be \emph{any} nonlinear continuous functions and $t$ (often the time) is
the independent variable. This kind of equations appear in many problems and in almost all scientific disciplines.

Equation (\ref{e01}) can be solved by the method of \emph{variation of parameters}: using $n$ linearly independent solutions, $y_1 (t), \cdots, y_n (t)$, of the homogenous part. Then, the general solution is just a linear combination of the independent solutions plus the particular solution associated to the nonhomogeneous equation \cite{Strang}. The variation of parameters method relies on the capability of finding the $n$ linearly independent solutions. Unfortunately, \emph{there is no general method} to find these solutions. Another method, called \emph{undetermined coefficients} \cite{Strang}, is restricted to the case of constant coefficient, only. Finally, in the specific case if $f (t)$ and all the $n + 1$ functions, $f_i (t)$, are polynomials, then an approximate solution can be found by power series \cite{Strang}. However, in this study, $f (t)$ and \emph{all} the $f_i (t)$ functions can be \emph{any} nonlinear continuous functions that are nonsingular in the integration time range.

The proposed Least-Squares (LS) method can be applied to solve Eq. (\ref{e01}) for \emph{any} value of $n$. However, without loss of generality and for sake of brevity, the approach is here applied to second order nonhomogeneous linear DE with nonconstant coefficients,
\begin{equation}\label{e02}
    f_2 (t) \, \dfrac{\dd^2 y (t)}{\dd t^2} + f_1 (t) \, \dfrac{\dd y (t)}{\dd t} + f_0 (t) \, y (t) = f (t).
\end{equation}

It is important to outline that, if functions $f_1 (t)$, $f_0 (t)$, and $f (t)$ are continuous and nonsingular within the integration range, the Initial Value Problems (IVP) always admit solutions, while Boundary Value Problems (BVP) may have a single, multiple, no, or infinite solutions. Final special analysis is dedicated to particular BVP (typical from optimal control) where the variable is a vector, $\{\B{x}\T, \; \B{\lambda}\T\}\T$, and where the state vector is defined at initial time, $\B{x} (t_0) = \B{x}_0$, and the costate vector at final time, $\B{\lambda} (t_f) = \B{\lambda}_f$.

\section{The constrained expressions}

The key idea of this study is to search the solution of Eq. (\ref{e01}) using \emph{constrained expressions}, whose theory is presented in Ref. \cite{Mortari}. These expressions have embedded \emph{all} the DE constraints, $\left.\dfrac{\dd^{\B{d}_i} y}{\dd t^{\B{d}_i}}\right|_{t = \B{t}_i} = y_{\B{t}_i}^{(\B{d}_i)}$, where the $n$-element vector, $\B{d}$, contains the constraints' derivatives orders and the $n$-element vector, $\B{t}$, indicates where the constraints are specified.

The \emph{constrained equations} adopted in this study are expressed as,
\begin{equation}\label{e03}
    \boxed{ y (t) = g (t) + \ds\sum_{i = 1}^n \beta_i (t, \B{t}) \left[y_{\B{t}_i}^{(\B{d}_i)} - g_{\B{t}_i}^{(\B{d}_i)}\right] \qquad \text{where:} \quad \beta_i^{(\B{d}_k)} (t_k, \B{t}) = \delta_{ik} }
\end{equation}
expressions that are \emph{linear} functions in the unknown function $g (t)$ and in its derivatives, $g_{\B{t}_i}^{(\B{d}_i)}$, evaluated at constraints times and where $\delta_{ik}$ is the Kronecker delta. The $\beta_i (t, \B{t})$ are special functions of the time and constraints times defined by the vector $\B{t}$. The $\beta_i$ functions given in Eq. (\ref{e03}) are not unique but they are characterized by $\beta_i^{(\B{d}_k)} (t_k, \B{t}) = \delta_{ik}$. Detailed derivations and presentations of these \emph{constrained expressions} can be found in Ref. \cite{Mortari}. However, let's give three constrained expression examples.

{\bf Example \#1.} In the first example consider the function,
\begin{equation}\label{e04}
    y (t) = g (t) + \dfrac{t \, (2 t_2 - t)}{2(t_2 - t_1)} \, (\dot{y}_1 - \dot{g}_1) + \dfrac{t \, (t - 2 t_1)}{2(t_2 - t_1)} \, (\dot{y}_2 - \dot{g}_2),
\end{equation}
where $\beta_1 (t, \B{t}) = \dfrac{t \, (2 t_2 - t)}{2(t_2 - t_1)}$ and $\beta_2 (t, \B{t}) = \dfrac{t \, (t - 2 t_1)}{2(t_2 - t_1)}$. The first derivative of Eq. (\ref{e04}) is
\begin{equation*}
    \dot{y} (t) = \dot{g} (t) + \dfrac{t_2 - t}{t_2 - t_1} (\dot{y}_1 - \dot{g}_1) + \dfrac{t - t_1}{t_2 - t_1} (\dot{y}_2 - \dot{g}_2).
\end{equation*}
It is easy to verify that, when $t = \B{t} (1) = t_1$ then $\dot{y} (t_1) = \dot{y}_1$ and when $t = \B{t} (2) = t_2$ then $\dot{y} (t_2) = \dot{y}_2$. Therefore, no matter what $g (t)$ is, Eq. (\ref{e04}) can be used as \emph{constrained expression} for functions subject to: $\dot{y} (t_1) = \dot{y}_1$ and $\dot{y} (t_2) = \dot{y}_2$. \\

{\bf Example \#2.} This example is for a function subject to the following $n = 4$ constraints,
\begin{equation*}
    \left.\dfrac{\dd^2 y}{\dd t^2}\right|_{\B{t}_1} = \ddot{y}_{\B{t}_1}, \qquad y (\B{t}_2) = y_{\B{t}_2}, \qquad y (\B{t}_3) = y_{\B{t}_3}, \qquad \text{and} \qquad \left.\dfrac{\dd y}{\dd t}\right|_{\B{t}_4} = \dot{y}_{\B{t}_4}.
\end{equation*}
where, $\B{d} = \{2, \; 0, \; 0, \; 1\}$. Let's select the constraint time vector as, $\B{t} = \{-1, \; 0, \; 2, \; 2\}$. A constrained expression with embedded all four constraints is
\begin{equation}\label{e05}
\begin{array}{ll}
    y (t) = g (t) & + \dfrac{- 4 + 4 \, t - t^2}{14} \, t \, (\ddot{y}_{\B{t}_1} - \ddot{g}_{\B{t}_1}) + \dfrac{28 - 24 \, t +  3 \, t^2 + t^3}{28} (y_{\B{t}_2} - g_{\B{t}_2}) + \\ ~ & + \dfrac{24 -  3 \, t - t^2}{28} \, t \, (y_{\B{t}_3} - g_{\B{t}_3}) + \dfrac{- 10 \, t +  3 \, t^2 + t^3}{14} (\dot{y}_{\B{t}_4} - \dot{g}_{\B{t}_4})
\end{array}
\end{equation}
where
\begin{equation*}
    \beta_1 (t, \B{t}) = \dfrac{- 4 + 4 t - t^2}{14} t, \qquad \beta_2 (t, \B{t}) = \dfrac{28 - 24 t +  3 t^2 + t^3}{28},
\end{equation*}
\begin{equation*}
    \beta_3 (t, \B{t}) = \dfrac{24 -  3 t - t^2}{28} t, \qquad \text{and} \qquad \beta_4 (t, \B{t}) = \dfrac{- 10 t +  3 t^2 + t^3}{14}.
\end{equation*}
It is not difficult to verify that $y (t)$, as defined by Eq. (\ref{e05}), has embedded all four constraints, $(\ddot{y}_{\B{t}_1}, \, y_{\B{t}_2}, \, y_{\B{t}_3}, \, \dot{y}_{\B{t}_4})$, independent what $g (t)$ is.

{\bf Example \#3.} This example shows the constrained expression when the constraints are specified in a relative way, as for
\begin{equation*}
    y (t_1) = y (t_2) \qquad \text{and} \qquad \dot{y} (t_1) = \dot{y} (t_2).
\end{equation*}
In this specific case, a constrained expression is
\begin{equation*}
    y (t) = g (t) + \dfrac{t}{t_2 - t_1} \, (g_1 - g_2) + \dfrac{t - 2 (t_1 + t_2)}{2(t_2 - t_1)} \, t \, (\dot{g}_1 - \dot{g}_2)
\end{equation*}
It is straightforward to prove this equation satisfies the two relative constraints, $y_1 = y_2$ and $\dot{y}_1 = \dot{y}_2$.

These three examples show that the solution, $y (t)$, can be expressed in term of an unknown function, $g (t)$, such that $y (t)$ always satisfies \emph{all} DE constraints. This allows us to re-write the original DE in terms of the new function, $g (t)$, thus obtaining a DE with constraints already embedded in the DE. This new DE has two interesting properties: \emph{it is not subject to external constraints} and \emph{it is linear in $g (t)$ and its derivatives}.

In this study simple constrained equations have been provided to solve linear DE with nonconstant coefficients for IVP, in Eqs. (\ref{e12}, \ref{e20}, \ref{e21}), and for BVP, in Eqs. (\ref{e23}, \ref{e29}, \ref{e30}, \ref{e31}, \ref{e32}, \ref{e33}, \ref{e34}, \ref{e35}, \ref{e36}).

\subsection{The least-squares approach}

Since function $g (t)$ is free to be selected, then it can be expressed as a linear combinations of a set of $m$ linearly independent basis functions, $h_k (t)$,
\begin{equation}\label{e06}
    g (t) = \B{\xi}\T \, \B{h} (t) = \ds\sum_{k = 0}^m \xi_k \, h_k (t).
\end{equation}
This means that two distinct functions, $h_i (t)$ and $h_j (t)$ with $i\ne j$, \emph{must span different function spaces}. In addition, all functions, $h_k (t)$, and their derivatives must be continuous and nonsingular within the time range. By doing this, the coefficients $\xi_k$ of Eq. (\ref{e06}) become our unknowns. Once these coefficients are computed, the $g (t)$ function is known and, consequently, Eq. (\ref{e04}) provides the DE solution.

Examples of basis functions are polynomials (e.g., Lagrange, Legendre, monomial, Chebyshev, etc.), Fourier series, monomial plus Fourier series, and any combinations of continuous and nonsingular functions spanning different function spaces.

By substituting the expression of $y (t)$, given in Eq. (\ref{e03}), in the DE of Eq. (\ref{e02}), along with the expression of $g (t)$, given in Eq. (\ref{e06}), and its derivatives, $\dot{g} (t) = \B{\xi}\T \dot{\B{h}} (t)$ and $\ddot{g} (t) = \B{\xi}\T \ddot{\B{h}} (t)$, a \emph{linear equation} in terms of the unknown coefficients, $\B{\xi}$, is obtained. This equation can be then specialized for a set of $N$ values of $t_j$ (e.g., uniformly distributed in the integration time range), obtaining
\begin{equation*}
    \ds\sum_{k = 0}^m \xi_k \, p_k (t_j) = \B{p}\T (t_j) \, \B{\xi} = \lambda (t_j)
\end{equation*}
where $\B{p} (t_j)$ is an $n$-long known vector and $j \in [1, N]$ and $N \ge m + 1$. This set of $N$ equations in $m + 1$ unknowns (usually, $N \gg m$) can be set in the matrix form
\begin{equation*}
    P \, \B{\xi} = \begin{bmatrix} p_0 (t_1) & p_1 (t_1) & \cdots & p_m (t_1)\\ p_0 (t_2) & p_1 (t_2) & \cdots & p_m (t_2)\\ \vdots & \vdots & \ddots & \vdots\\ p_0 (t_N) & p_1 (t_N) & \cdots & p_m (t_N)\end{bmatrix} \, \begin{Bmatrix} \xi_0\\ \xi_1\\ \vdots\\ \xi_m\end{Bmatrix} = \begin{Bmatrix} \lambda (t_1)\\ \lambda (t_2)\\ \vdots\\ \lambda (t_N)\end{Bmatrix} = \B{\lambda}
\end{equation*}
admitting the LS solution
\begin{equation}\label{e07}
    \boxed{ \B{\xi} = (P\T P)^{-1} \, P\T \, \B{\lambda} }
\end{equation}
This LS solution is computed by scaling the $P$ matrix is order to decrease the condition number of $P\T P$ and, consequently, the numerical errors. This procedure is applied in detail for a second order nonhomogeneous linear DE with nonconstant coefficients for IVP and BVP, respectively.

\subsection{The selected basis functions}

\begin{figure}[ht]
    \centering\includegraphics[width=\linewidth]{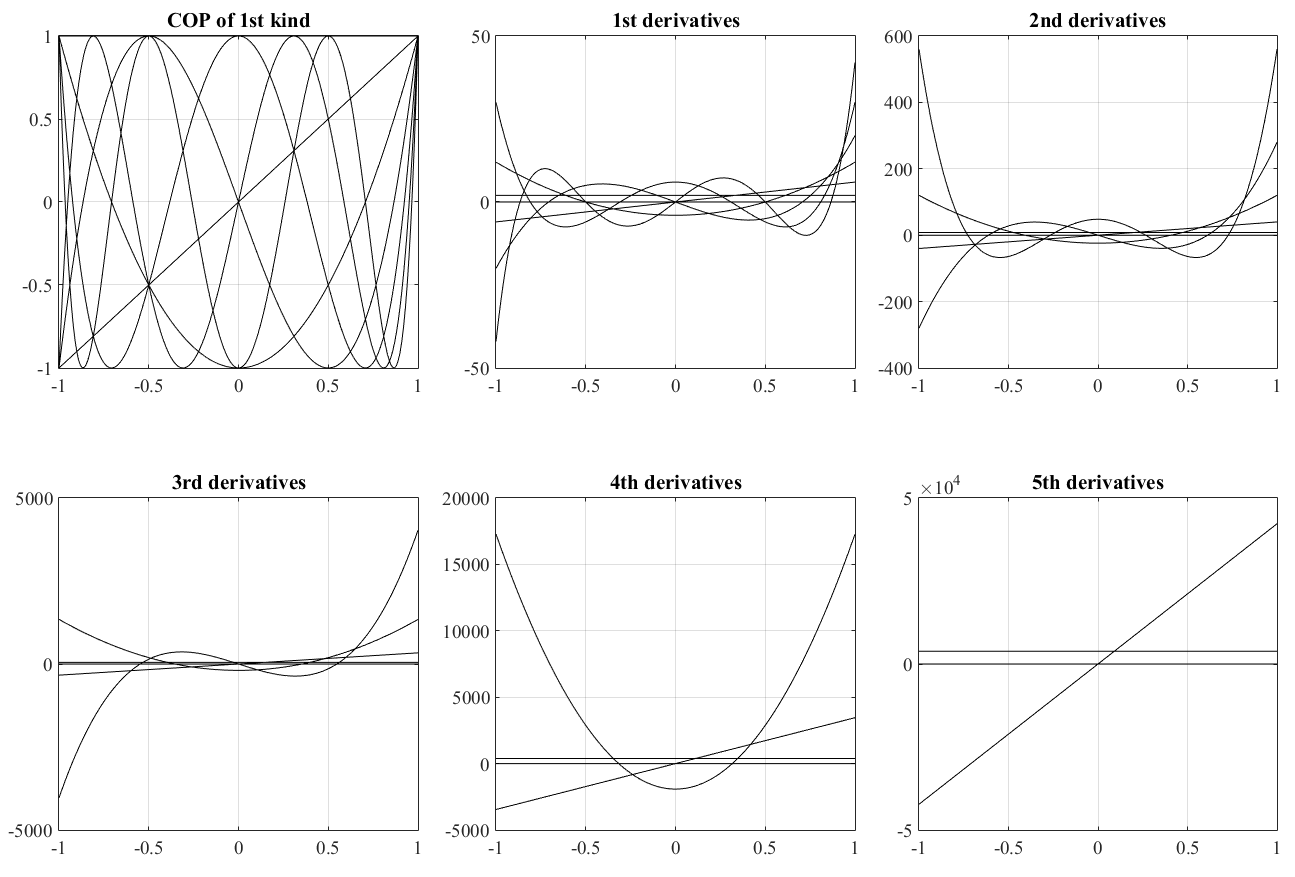}
	\caption{Chebyshev Orthogonal Polynomials and the first two derivatives}
    \label{fig00}
\end{figure}

In all the examples provided in this paper, the Chebyshev Orthogonal Polynomials (COP) of the first kind have been selected to represent the basis functions set. Note that, this selection \emph{may not be the best solution}. In fact, while COP are a versatile basis to describe almost any kind of function, the COP derivatives are affected by a sort of Runge's phenomenon, with one order of magnitude increase at each subsequent derivatives. Figure \ref{fig00} shows this effect for the first two derivatives.

Since COP are defined in term of a new variable, $x \in [-1, +1]$, we set $x$ linearly related to $t \in [t_1, t_2]$, as
\begin{equation}\label{e08}
    x = 2 \dfrac{t - t_1}{t_2 - t_1} - 1 \qquad \longleftrightarrow \qquad t = t_1 + \dfrac{(x + 1)(t_2 - t_1)}{2},
\end{equation}
where $t_2$ is specifically defined in BVP while it can be considered as the integration upper limit in IVP. Setting $\delta t = t_2 - t_1$, the derivatives in terms of the new variable are,
\begin{equation}\label{e09}
    \begin{array}{ll}
        \dfrac{\dd y}{\dd t} =& \dfrac{\dd y}{\dd x} \cdot \dfrac{\dd x}{\dd t} = \dfrac{2}{\delta t} \, \dfrac{\dd y}{\dd x} \\
        \dfrac{\dd^2 y}{\dd t^2} =& \dfrac{\dd}{\dd t} \left(\dfrac{\dd y}{\dd t}\right) = \dfrac{\dd}{\dd x} \left(\dfrac{\dd y}{\dd t}\right) \cdot \dfrac{\dd x}{\dd t} = \dfrac{\dd}{\dd x} \left(\dfrac{2}{\delta t} \, \dfrac{\dd y}{\dd x}\right) \cdot \dfrac{\dd x}{\dd t} = \dfrac{4}{\delta t^2} \, \dfrac{\dd^2 y}{\dd x^2}
    \end{array}.
\end{equation}
Therefore, Eq. (\ref{e02}) can be re-written as,
\begin{equation}\label{e10}
    \dfrac{4}{\delta t^2} f_2 (x) \, \dfrac{\dd^2 y (x)}{\dd x^2} + \dfrac{2}{\delta t} \, f_1 (x) \, \dfrac{\dd y (x)}{\dd x} + f_0 (x) \, y (x) = f (x),
\end{equation}
where the functions $f_2$, $f_1$, $f_0$, and $f$ are now expressed in term of the new variable using Eq. (\ref{e08}). By changing the integration variable, particular attention must be given to the constraints specified in term of derivatives. In fact, the derivatives $\dfrac{\dd y}{\dd t}$ and $\dfrac{\dd^2 y}{\dd t^2}$ are related to derivatives $\dfrac{\dd y}{\dd x}$ and $\dfrac{\dd^2 y}{\dd x^2}$ as specified by Eq. (\ref{e09}). Therefore, also the constraints, provided in term of the first and/or second derivatives, need to comply with the rules given in Eq. (\ref{e09}),
\begin{equation*}
    \left.\dfrac{\dd y}{\dd x}\right|_{x_1} = \left.\dfrac{\dd y}{\dd t}\right|_{t_1} \dfrac{\delta t}{2} = \dfrac{\delta t}{2} \, \dot{y}_1 = \dot{y}_{1x} \quad \text{and} \quad \left.\dfrac{\dd^2 y}{\dd x^2}\right|_{x_1} = \left.\dfrac{\dd^2 y}{\dd t^2}\right|_{t_1} \dfrac{\delta t^2}{4} = \dfrac{\delta t^2}{4} \, \ddot{y}_1 = \ddot{y}_{1x},
\end{equation*}
meaning that: \emph{the constraints on the derivatives in term of the new $x$ variable now depend on the integration time range}.

\section{Least-squares Solution of Initial Value Problems}

Three distinct IVPs can be considered, depending on the DE constraints kind. Consider first the most classic problem where the function and its first derivative are specified in one point.

\subsection{Initial Value Problems subject to: $y(t_1) = y_1$ and $\dot{y} (t_1) = \dot{y}_1$}

In this case, the constraints written in term of the new variable ($x$) are,
\begin{equation*}
    y(x_1 =-1) = y_1 \qquad \text{and} \qquad \left.\dfrac{\dd y}{\dd x}\right|_{x_1 =-1} = \dot{y}_{1x} = \dfrac{\delta t}{2} \, \dot{y}_1,
\end{equation*}
and a simple \emph{constrained expression} for this IVP is,
\begin{equation}\label{e11}
    y (x) = g(x) + (y_1 - g_1) + (x + 1) (\dot{y}_{1x} - \dot{g}_1),
\end{equation}
where $g (-1) = g_1$ and $\dot{g} (-1) = \dot{g}_1$. Again, the solution $y (x)$, as expressed by Eq. (\ref{e11}), has embedded the constraints, no matter what the function $g (x)$ is.  Substituting $y (x)$, as expressed by Eq. (\ref{e11}), in Eq. (\ref{e10}) we obtain,
\begin{equation}\label{e12}
    \dfrac{4}{\delta t^2} f_2 \dfrac{\dd^2 g}{\dd x^2} + \dfrac{2}{\delta t} f_1 \left(\dfrac{\dd g}{\dd x} - \dot{g}_1\right) + f_0 \left[g - g_1 - \dot{g}_1 (x + 1)\right] = f - \dfrac{2}{\delta t} \dot{y}_{1x} f_1 - f_0 \left[y_1 + \dot{y}_{1x} (x + 1)\right].
\end{equation}
Now, let $g (x)$ be expressed as a linear combination of COPs of the first kind,
\begin{equation}\label{e13}
    g (x) = \ds\sum_{k = 0}^m \xi_k \, T_k (x),
\end{equation}
which are defined by the recursive function,
\begin{equation}\label{e14}
    T_{k + 1} = 2 \, x \, T_k - T_{k - 1} \qquad \text{starting from:} \; \left\{\begin{array}{ll} T_0 =& 1 \\ T_1 =& x\end{array}\right. .
\end{equation}
All derivatives of COP can be computed in a recursive way, starting from
\begin{equation*}
    \dfrac{\dd T_0}{\dd x} = 0, \quad \dfrac{\dd T_1}{\dd x} = 1 \qquad \text{and} \qquad \dfrac{\dd^d T_0}{\dd x^d} = \dfrac{\dd^d T_1}{\dd x^d} = 0 \quad (\forall \; d > 1),
\end{equation*}
while the subsequent derivatives of Eq. (\ref{e14}) give for $k > 1$,
\begin{equation}\label{e15}
    \begin{array}{cccccc}
        \dfrac{\dd T_{k+1}}{\dd x} &=& 2 \, T_k &+ 2 x \, \dfrac{\dd T_k}{\dd x} &- \dfrac{\dd T_{k-1}}{\dd x} \\ [8pt]
        \dfrac{\dd^2 T_{k+1}}{\dd x^2} &=& 4 \dfrac{\dd T_k}{\dd x} &+ 2 x \, \dfrac{\dd^2 T_k}{\dd x^2} &- \dfrac{\dd^2 T_{k-1}}{\dd x^2} \\ [4pt]
        \vdots & ~ & \vdots & \vdots & \vdots \\ [4pt]
        \dfrac{\dd^d T_{k+1}}{\dd x^d} &=& 2 d \, \dfrac{\dd^{d-1} T_k}{\dd x^{d-1}} &+ 2 x \, \dfrac{\dd^d T_k}{\dd x^d} &- \dfrac{\dd^d T_{k-1}}{\dd x^d}, & (\forall \; d \ge 1),
    \end{array}
\end{equation}
In particular, it is easy to show that,
\begin{equation}\label{e16}
    T_k (-1) = (-1)^k, \quad \left.\dfrac{\dd T_k}{\dd x}\right|_{x=-1} = (-1)^{k+1} \, k^2, \quad \left.\dfrac{\dd^2 T_k}{\dd x^2}\right|_{x=-1} = (-1)^k \, \dfrac{k^2 \, (k^2 - 1)}{3}.
\end{equation}
Therefore, substituting the expressions given in Eqs. (\ref{e13}-\ref{e16}) in Eq. (\ref{e12}), the following equation
\begin{equation}\label{e17}
    \begin{array}{c} \ds\sum_{k = 0}^m \xi_k \left\{\dfrac{4}{\delta t^2} f_2 \dfrac{\dd^2 T_k}{\dd x^2} + \dfrac{2}{\delta t} f_1 \left[\dfrac{\dd T_k}{\dd x} - (-1)^{k+1} k^2\right] + f_0 \left[T_k - (-1)^k - (-1)^{k+1} k^2 (x + 1)\right]\right\} = \\ = f - \dfrac{2}{\delta t} \dot{y}_{1x} f_1 - f_0 [y_1 + \dot{y}_{1x} (x + 1)] \} \end{array}
\end{equation}
is obtained. However, particular attention must be given to Eq. (\ref{e17}) because, for $k = 0$ and $k = 1$, all three terms of the RHS vanish,
\begin{equation*}
    \dfrac{\dd^2 T_k}{\dd x^2} = \dfrac{\dd T_k}{\dd x} - (-1)^{k+1} k^2 = T_k - (-1)^k - (-1)^{k+1} k^2 (x + 1) = 0 \qquad \text{for } k = 0 \;\; \text{and} \;\; k = 1
\end{equation*}
which is equivalent to rewrite Eq. (\ref{e17}) as
\begin{equation}\label{e18}
    \begin{array}{c} \ds\sum_{k = 2}^m \xi_k \left\{\dfrac{4}{\delta t^2} f_2 \dfrac{\dd^2 T_k}{\dd x^2} + \dfrac{2}{\delta t} f_1 \left[\dfrac{\dd T_k}{\dd x} - (-1)^{k+1} k^2\right] + f_0 \left[T_k - (-1)^k - (-1)^{k+1} k^2 (x + 1)\right]\right\} = \\ = f - \dfrac{2}{\delta t} \dot{y}_{1x} f_1 - f_0 [y_1 + \dot{y}_{1x} (x + 1)] \}\end{array}
\end{equation}
The reason why in Eq. (\ref{e17}) for $k = 0$ and $k = 1$, all three terms of the RHS vanish, derives from the fact that the first two terms of COP are constant and linear in $x$. Now, the constrained expression of Eq. (\ref{e11}) is derived using a constant plus a linear expression in $x$. This means that the basis functions used for $g (x)$ cannot be composed using the same function spaces, namely, the constant and the linear expression in $x$, because already adopted to define the constrained expression.

Note that, the two derivatives, $\dfrac{\dd^2 T_k}{\dd x^2}$ and $\dfrac{\dd T_k}{\dd x}$, are specific known polynomials in $x$, derived using Eqs. (\ref{e13}-\ref{e15}). Therefore, Eq. (\ref{e18}) is a linear equation in terms of the $(m - 1)$ unknown coefficients $\xi_k$. This allows us to estimate these coefficients by LS, by specifying Eq. (\ref{e18}) for a set of $N$ values $x_j$, ranging from $x_1 =-1$ to $x_2 =+1$. Specifically, we have
\begin{align*}
    p_k (x_j) =& \dfrac{4}{\delta t^2} f_2 (x_j) \left.\dfrac{\dd^2 T_k}{\dd x^2}\right|_{x_j} + \dfrac{2}{\delta t} f_1 (x_j) \left[\left.\dfrac{\dd T_k}{\dd x}\right|_{x_j} - (-1)^{k+1} k^2\right] + \\ ~& + f_0 (x_j) \left[T_k (x_j) - (-1)^k - (-1)^{k+1} k^2 (x_j + 1)\right] \\
    \lambda (x_j) =& f (x_j) - \dfrac{2}{\delta t} \dot{y}_{1x} f_1 (x_j) - f_0 (x_j) [y_1 + \dot{y}_{1x} (x_j + 1)]
\end{align*}

In the next subsection the proposed approach is applied to a DE with known analytical solution. Accuracy comparison is provided with respect to the solution obtained by the Runge-Kutta-Fehlberg step-varying integrator (MATLAB function \verb"ODE45").

\subsection{Accuracy tests}

\begin{figure}[ht]
    \centering\includegraphics[width=\linewidth]{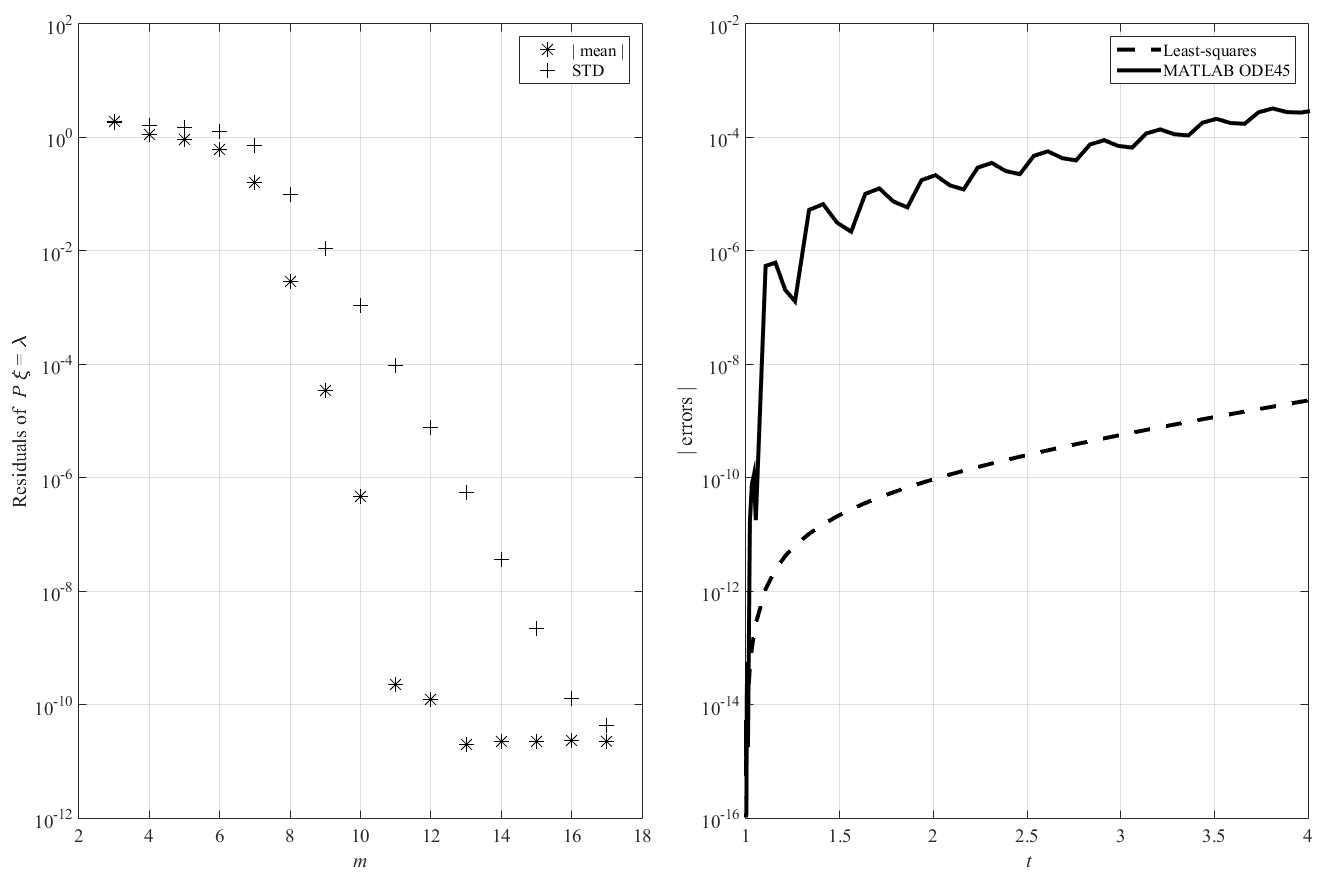}
	\caption{Results example with known solution ($m = 16$ and $N = 1,000$)}
    \label{fig1}
\end{figure}

Consider integrating the following IVP from $t_1 = 1$ to $t_2 = 4$,
\begin{equation}\label{e19}
    t^2 \, \dfrac{\dd^2 y}{\dd t^2} - t(t + 2) \, \dfrac{\dd y}{\dd t} + (t + 2) \, y = 0 \qquad \text{subject to:} \; \left\{\begin{array}{l} y (1) = y_1 = 1\\ \dot{y} (1) = \dot{y}_1 = 0\end{array}\right.
\end{equation}
This implies $\dot{y}_{1x} = \dot{y}_1 \dfrac{3}{2} = 0$. The general solution of this equation is $y (t) = \left(2 - e^{t - 1}\right) \, t$. Equation (\ref{e19}) has been solved using the proposed LS approach (with $m = 16$ and $N = 1,000$) and integrated using MATLAB function \verb"ODE45", implementing the Runge-Kutta-Fehlberg variable step integrator. The results are shown in Fig. \ref{fig1}.

In the left plot of Fig. \ref{fig1} the absolute values of mean and standard deviation of the $(P \, \B{\xi} - \B{\lambda})$ residuals are shown as a function of $m$. When the residuals standard deviation reaches the minimum (at $m = 17$)\footnote{The value of $m = 17$ implies a $16\times 16$ size of matrix $(P\T P)$.} the LS approach provides the best accuracy results. The errors with respect to the true solution for the LS approach and the errors obtained using MATLAB function \verb"ODE45", are shown in the right plot. For this IVP, the LS method provides about \emph{five order of magnitudes accuracy gain} with respect to \verb"ODE45" integrator.

\subsection{Initial Value Problems subject to: $y (t_1) = y_1$ and $\ddot{y} (t_1) = \ddot{y}_1 \; \to \; \ddot{y}_{1x} = \ddot{y}_1 \, \dfrac{\delta t^2}{4}$}

Using the constrained equation,
\begin{equation*}
    y (x) = g (x) - x \, (y_1 - g_1) + \dfrac{x^2 + x}{2} (\ddot{y}_{1x} - \ddot{g}_1)
\end{equation*}
then Eq. (\ref{e10}) becomes
\begin{equation}\label{e20}
\begin{array}{c}
     \dfrac{4}{\delta t^2} f_2 \left(\dfrac{\dd^2 g}{\dd x^2} - \ddot{g}_1\right) + \dfrac{2}{\delta t} f_1 \left(\dfrac{\dd g}{\dd x} + g_1 - \ddot{g}_1 \dfrac{2 x + 1}{2}\right) + f_0 \left(g + g_1 x - \ddot{g}_1 \dfrac{x^2 + x}{2}\right) = \\ = f - \dfrac{4}{\delta t^2} f_2 \, \ddot{y}_{1x} - \dfrac{2}{\delta t} f_1 \left(-y_1 + \ddot{y}_{1x} \dfrac{2 x + 1}{2}\right) - f_0 \left(- y_1 x + \ddot{y}_{1x} \dfrac{x^2 + x}{2}\right)
\end{array}
\end{equation}
Note that, if $y_1$ and $\ddot{y}_{1x}$ are known, then $\dot{y}_{1x}$ can be derived using Eq. (\ref{e10}) evaluated at $x_1 =-1$
\begin{equation*}
    \dot{y}_{1x} = \dfrac{\delta t^2 \, f (t_1) - 4 f_2 \, \ddot{y}_{1x} - \delta t^2 \, f_0 (t_1) \, y_1}{2 \, \delta t \, f_1 (t_1)}
\end{equation*}
provided that $f_1 (t_1) \ne 0$. Therefore, the solution of the DE given in Eq. (\ref{e10}) with constraints $y_1$ and $\ddot{y}_{1x}$ can be solved as in the previous section with constraints $y_1$ and $\dot{y}_{1x}$, where $\dot{y}_{1x}$ is provided as a function of $\dot{y}_1$ and integration time range $\delta t$.

\subsection{Initial Value Problems subject to: $\dot{y} (t_1) = \dot{y}_1$ and $\ddot{y} (t_1) = \ddot{y}_1$}

Using the constrained equation,
\begin{equation*}
    y (x) = g (x) + x \, (\dot{y}_{1x} - \dot{g}_1) + \left(\dfrac{x^2}{2} + x\right) (\ddot{y}_{1x} - \ddot{g}_1)
\end{equation*}
Eq. (\ref{e10}) becomes
\begin{equation}\label{e21}
\begin{array}{c}
    \dfrac{4}{\delta t^2} f_2 \left(\dfrac{\dd^2 g}{\dd x^2} - \ddot{g}_1\right) + \dfrac{2}{\delta t} f_1 \left[\dfrac{\dd g}{\dd x} - \dot{g}_1 - \ddot{g}_1 (x + 1)\right] + f_0 \left[g - \dot{g}_1 x - \ddot{g}_1 \left(\dfrac{x^2}{2} + x\right)\right] = \\ = f - \dfrac{4}{\delta t^2} f_2 \, \ddot{y}_{1x} - \dfrac{2}{\delta t} f_1 \left[\dot{y}_{1x} + \ddot{y}_{1x} (x + 1)\right] - f_0 \left[\dot{y}_{1x} \, x + \ddot{y}_{1x} \left(\dfrac{x^2}{2} + x\right)\right]
\end{array}
\end{equation}
Again, if $\dot{y}_{1x}$ and $\ddot{y}_{1x}$ are known, then $y_1$ can also be computed by specializing Eq. (\ref{e10}) at $x_1 =-1$
\begin{equation}\label{e22}
    y_1 = \dfrac{\delta t^2 \, f (t_1) - 4 f_2 \, \ddot{y}_{1x} - 2 \, \delta t \, f_1 (t_1) \, \dot{y}_{1x}}{\delta t^2 \, f_0 (t_1)} \qquad f_0 (t_1) \ne 0.
\end{equation}
Therefore, the solution of the DE given in Eq. (\ref{e10}) with constraints $\dot{y}_{1x}$ and $\ddot{y}_{1x}$ can be solved as in the previous section with constraints $y_1$ and $\dot{y}_{1x}$, where $y_1$ is provided by Eq. (\ref{e22}).

\section{Least-squares Solution of Boundary Value Problems}

Boundary Value Problems (BVP) appear in many applications arising in science and engineering. Examples are the modeling of chemical reactions, heat transfer, and diffusion. A thorough survey of the existing solutions to this problem can be found in Ref. \cite{Lin}, describing most of the existing methods with the exception of those using B\'ezier curves. The use of implicit B\'ezier functions to obtain approximate solutions of BVP (or IVP) is not a new idea, albeit it is quite recent (2004). Venkataraman has attacked the problem using optimization techniques \cite{Venkataraman1, Venkataraman2, Venkataraman3} while Zheng uses analytical LS approach \cite{Zheng}. B\'ezier curves have been adopted also to solve specific problems such as singular perturbed BVP \cite{Evrenosoglu} as well as integro-DE \cite{Ghomanjani}. Two-point BVP are usually solved by iterative techniques. The most common approaches are the shooting methods, transforming the BVP into IVP.

With respect to the analytical LS approach proposed in Ref. \cite{Zheng}, this paper has developed a practical, fast and easy to implement, numerical LS approach. The proposed method does not require any sophisticated optimization technique to solve BVP applied to linear second-order nonhomogeneous DE. Examples are provided with particularly emphasis on the approximate solutions accuracy levels.

Let's first consider the most common BVP whose constraints are $y (-1) = y_1$ and $y (1) = y_2$. A constrained function is
\begin{equation*}
    y (x) = g(x) + \dfrac{1 - x}{2} \, (y_1 - g_1) + \dfrac{1 + x}{2} \, (y_2 - g_2)
\end{equation*}
Substituting in Eq. (\ref{e10}) we obtain
\begin{equation}\label{e23}
    \begin{array}{c}
        \dfrac{4}{\delta t^2} f_2 \dfrac{\dd^2 g}{\dd x^2} + \dfrac{2}{\delta t} f_1 \left(\dfrac{\dd g}{\dd x} + \dfrac{g_1 - g_2}{2}\right) + f_0 \left(g - \dfrac{g_1 + g_2}{2} + \dfrac{g_1 - g_2}{2} x\right) = \\ = f - \dfrac{1}{\delta t} f_1 (y_2 - y_1) - f_0 \left(\dfrac{y_1 + y_2}{2} - \dfrac{y_1 - y_2}{2} x\right)
    \end{array}
\end{equation}
In particular, using COP to describe $g (x)$, we have
\begin{equation*}
    T_k (1) = 1, \qquad \left.\dfrac{\dd T_k}{\dd x}\right|_{x = 1} = k^2, \qquad \text{and} \qquad \left.\dfrac{\dd^2 T_k}{\dd x^2}\right|_{x = 1} = \dfrac{k^2 \, (k^2 - 1)}{3}.
\end{equation*}
Using these expressions in Eq. (\ref{e23}) we obtain
\begin{equation}\label{e24}
    \begin{array}{c}
        \ds\sum_{k = 0}^n \xi_k \left\{\dfrac{4}{\delta t^2} f_2 \dfrac{\dd^2 T_k}{\dd x^2} + \dfrac{2}{\delta t} f_1 \left[\dfrac{\dd T_k}{\dd x} + \dfrac{(-1)^k - 1}{2}\right] + f_0 \left[T_k - \dfrac{(-1)^k + 1}{2} + \dfrac{(-1)^k - 1}{2} x\right]\right\} = \\ = f - \dfrac{1}{\delta t} f_1 (y_2 - y_1) - f_0 \left(\dfrac{y_1 + y_2}{2} - \dfrac{y_1 - y_2}{2} x\right)
    \end{array}
\end{equation}
However, for $k = 0$ and $k = 1$, all three terms on the RHS of Eq. (\ref{e24}) becomes zeros
\begin{equation*}
    \dfrac{\dd^2 T_k}{\dd x^2} = \dfrac{\dd T_k}{\dd x} + \dfrac{(-1)^k - 1}{2} = T_k - \dfrac{(-1)^k + 1}{2} + \dfrac{(-1)^k - 1}{2} x = 0 \qquad \text{for} \; k = 0, 1
\end{equation*}
For this reason Eq. (\ref{e24}) must be selected as
\begin{equation}\label{e25}
    \begin{array}{c}
        \ds\sum_{k = 2}^m \xi_k \left\{\dfrac{4}{\delta t^2} f_2 \dfrac{\dd^2 T_k}{\dd x^2} + \dfrac{2}{\delta t} f_1 \left[\dfrac{\dd T_k}{\dd x} + \dfrac{(-1)^k - 1}{2}\right] + f_0 \left[T_k - \dfrac{(-1)^k + 1}{2} + \dfrac{(-1)^k - 1}{2} x\right]\right\} = \\ = f - \dfrac{1}{\delta t} f_1 (y_2 - y_1) - f_0 \left(\dfrac{y_1 + y_2}{2} - \dfrac{y_1 - y_2}{2} x\right)
    \end{array}
\end{equation}
The $(m-1)$ coefficients $\xi_k$ of Eq. (\ref{e25}) are then computed by LS using Eq. (\ref{e07}).

\subsection{Numerical accuracy tests with known solution}

Consider the following BVP with constant coefficients,
\begin{equation}\label{e26}
    \dfrac{\dd^2 y}{\dd t^2} + 2 \dfrac{\dd y}{\dd t} + y = 0 \qquad \text{subject to:} \; \left\{\begin{array}{l} y (0) = 1\\ y (1) = 3\end{array}\right.,
\end{equation}
whose solution is, $y (t) = e^{-t} + (3 e - 1) t e^{-t}$, with derivatives, $\dot{y} (t) =-e^{-t} (3 e t - t - 3 e + 2)$, and $\ddot{y} (t) = e^{-t} (3 e t - t - 6 e + 3)$.

\begin{figure}[ht]
    \centering\includegraphics[width=\linewidth]{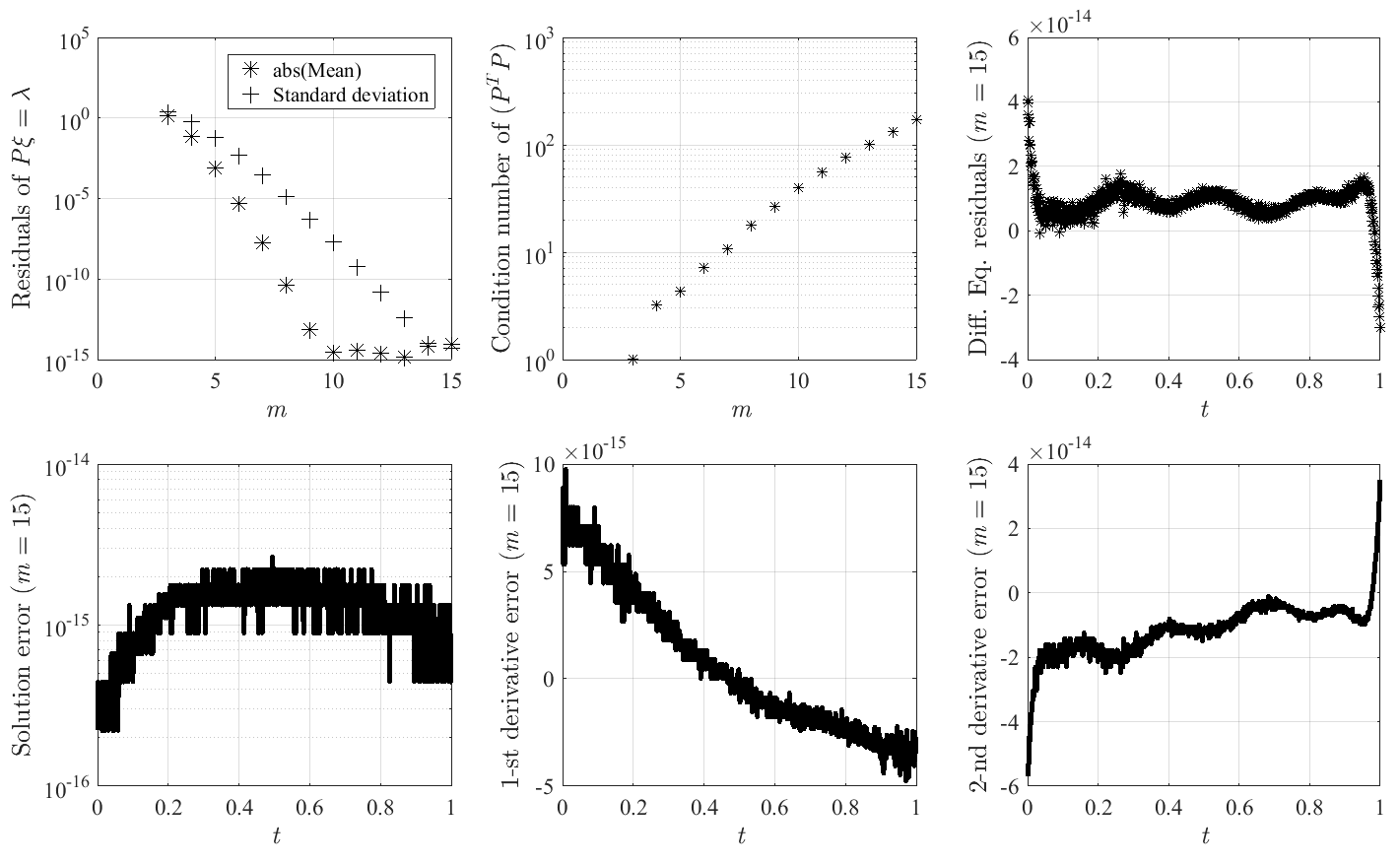}
	\caption{IVP least-square results for Eq. (\ref{e26}) and for $m\in[3,23]$}
    \label{fig2}
\end{figure}

Figure \ref{fig2} show the LS approach results for this test in terms of mean and standard deviation of $(P \, \B{\xi} - \B{\lambda})$ residuals (top, left) and the condition number of matrix $P\T P$ (top, center) as a function of the number of COP adopted ($m-1$) to solve the LS problem. The residuals of the DE of Eq. (\ref{e26}) are provided (top, right), and the errors of the LS solution, $y (t)$ (bottom, left), the first derivative, $\dot{y} (t)$ (bottom, center), and the second derivative, $\ddot{y} (t)$ (bottom, right), with respect to the true solution, are provided.

Shooting methods are transforming BVP into IVP. Numerical integrations of IVP, provide subsequent estimates based on previous estimates. This implies that the error, in general, is accumulating. In the contrary, the accuracy provided by LS solution is ``uniformly'' distributed within the integration bounds. If more accuracy is desired on a specific range, then by increasing the number of points on that range or by providing greater weights to the points on that range, the accuracy increase is obtained where desired.

\subsection{Tests with unknown solution, with no solution, and with infinite solutions}

\begin{figure}[ht]
    \centering\includegraphics[width=\linewidth]{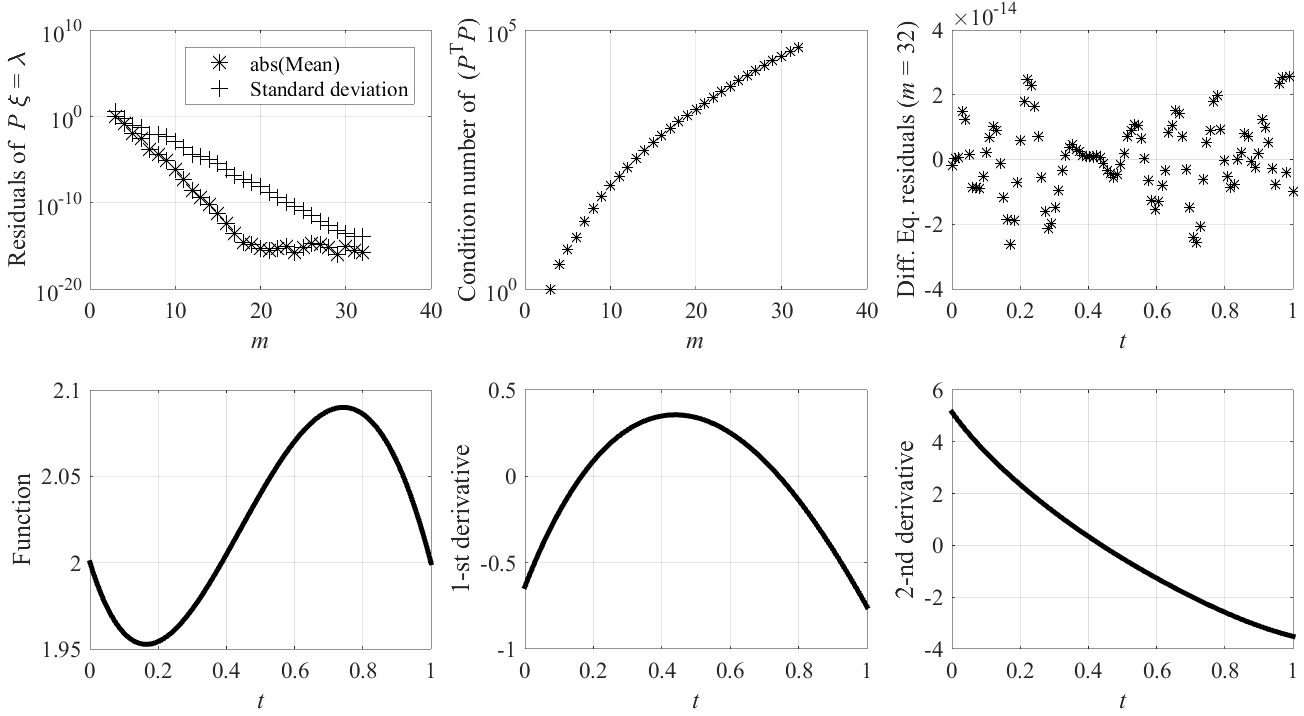}
	\caption{Results with unknown solution}
    \label{fig3}
\end{figure}

Consider the DE with \emph{unknown solution},
\begin{equation*}
    (1 + 2t) \dfrac{\dd^2 y}{\dd t^2} + \left(\cos t^2 - 3 t + 1\right) \dfrac{\dd y}{\dd t} + \left(6 \sin t^2 - e^{\cos (3 t)}\right) y = \dfrac{2[1 - \sin(3 t)](3 t - \pi)}{4 - t},
\end{equation*}
subject to $y (0) = 2$ and $y (1) = 2$. In this case the LS solution results are given in the plots of Fig. \ref{fig3} with the same meaning to those provided in Fig. \ref{fig2}. The LS solution accuracy increases up to $m = 21$-degree COP. The standard deviation of the residuals reaches about $10^{-14}$ accuracy level while the DE residuals are lower than $4.0 \cdot 10^{-14}$.

Consider the following BVP with \emph{no solution},
\begin{equation}\label{e27}
    \dfrac{\dd^2 y}{\dd t^2} - 6 \, \dfrac{\dd y}{\dd t} + 25 \, y = 0 \qquad \text{subject to:} \; \left\{\begin{array}{l} y (0) = 1\\ y (\pi) = 2\end{array}\right..
\end{equation}
In fact, the general solution of Eq. (\ref{e27}) is $y (t) = [a \, \cos(4\, t) + b \, \sin(4\, t)] \, e^{3\, t}$, where the constraint $y (0) = 1$ gives $a = 1$ while the constraint $y (\pi) = 2$ gives $2 = e^{3 \, \pi}$, a wrong identity!
\begin{figure}[ht]
    \centering\includegraphics[width=\linewidth]{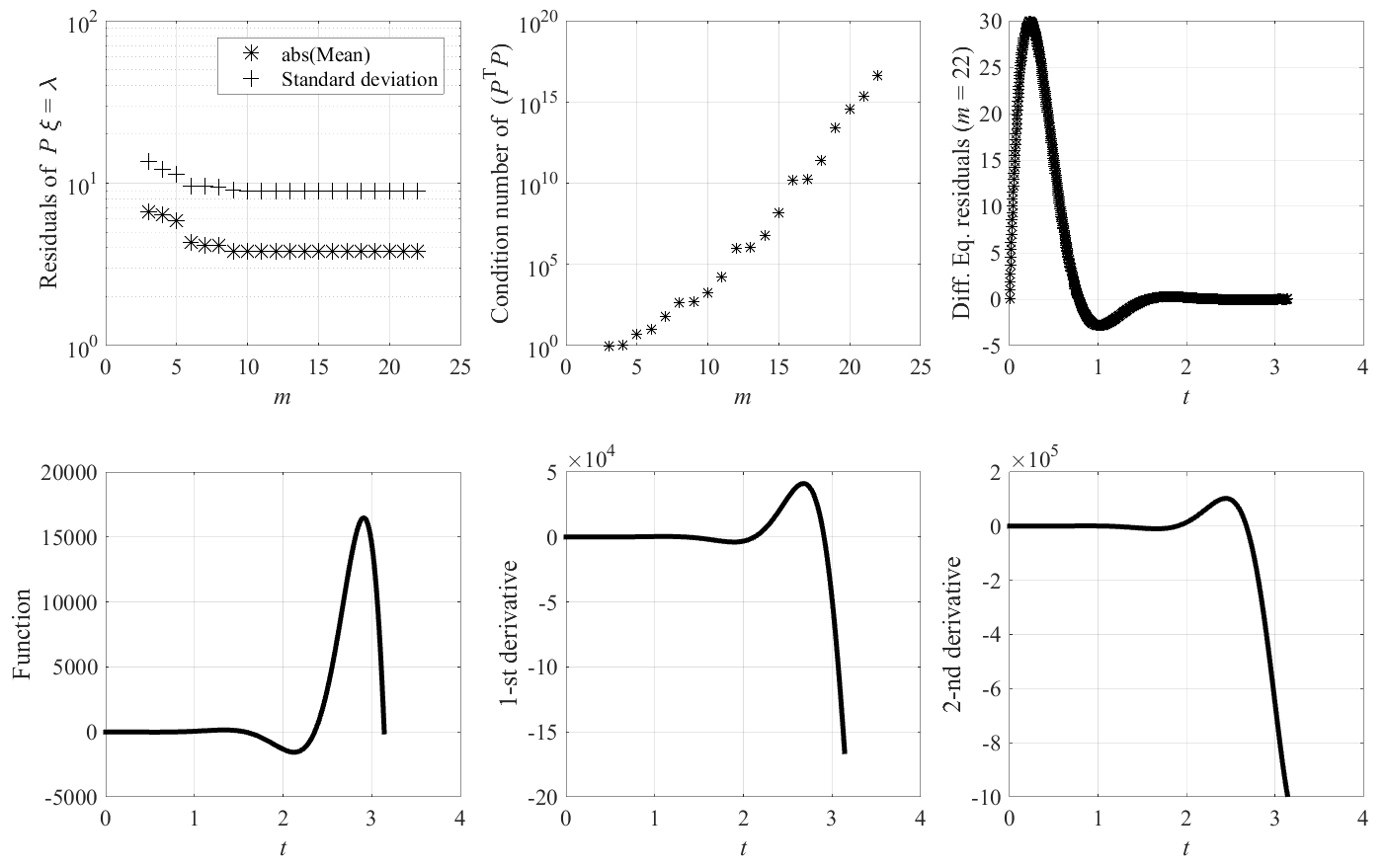}
	\caption{Results with NO solution}
    \label{fig4}
\end{figure}

Figure \ref{fig4} shows the results when trying to solve by LS the problem given in Eq. (\ref{e27}). For this example the number of COP terms has been increased up to when the matrix $P\T P$ become numerically singular, at $m = 22$, with a condition number value above $10^{15}$. The mean and standard deviation of the $(P \, \B{\xi} - \B{\lambda})$ residuals show no convergence while the condition number of $P\T P$ indicates that the problem has no solution. It is possible to show that, even in the no solution case, the proposed LS approach provides anyway a ``solution'' complying with the DE constraints!

\begin{figure}[ht]
    \centering\includegraphics[width=\linewidth]{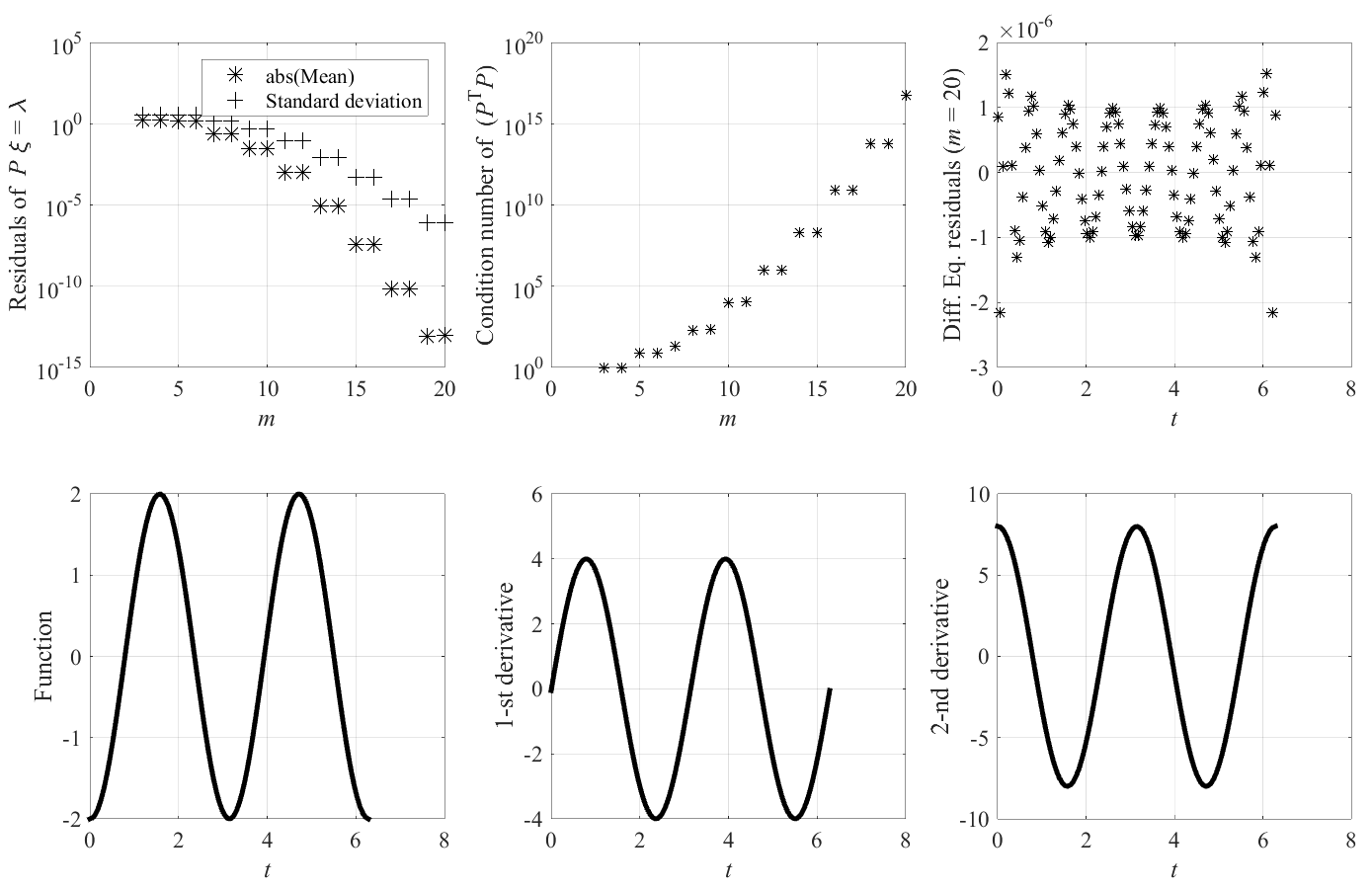}
	\caption{Results with infinite solutions}
    \label{fig5}
\end{figure}

Finally, consider the BVP with \emph{infinite solutions},
\begin{equation}\label{e28}
    \dfrac{\dd^2 y}{\dd t^2} + 4 \, y = 0 \qquad \text{subject to:} \; \left\{\begin{array}{l} y (0) =-2\\ y (2\pi) =-2\end{array}\right..
\end{equation}
In fact, the general solution of Eq. (\ref{e28}) is, $y (t) = a \, \cos(2 \, t) + b \, \sin (2 \, t)$, consisting of infinite solutions as $b$ can have any value. Results of the LS approach are given in Fig. \ref{fig5} showing the convergence, but not at machine level accuracy.
Note the differences between the two cases of no and infinite solutions. Both of them experience bad condition number, but the convergence is experienced in the infinite solution case, only.

\subsection{Constraints: $y (t_1) = y_1$ and $\dot{y} (t_2) = \dot{y}_2 \; \to \; \dot{y}_{2x} = \dot{y}_2 \, \dfrac{2}{\delta t}$}

For this case the constrained equation
\begin{equation*}
    y (x) = g (x) + (y_1 - g_1) + (x + 1) (\dot{y}_{2x} - \dot{g}_2)
\end{equation*}
can be used. Substituting this equation in Eq. (\ref{e10}), we obtain
\begin{equation}\label{e29}
    \begin{array}{c} \dfrac{4}{\delta t^2} f_2 \dfrac{\dd^2 g}{\dd x^2} + \dfrac{2}{\delta t} f_1 \left(\dfrac{\dd g}{\dd x} - \dot{g}_2\right) + f_0 \left[g - g_1 - \dot{g}_2 (x + 1)\right] = \\ = f - \dfrac{2}{\delta t} f_1 \dot{y}_{2x} - f_0 \left[y_1 + \dot{y}_{2x} (x + 1)\right]\end{array}
\end{equation}
Then, using the expressions provided in Eqs. (\ref{e13}-\ref{e16}) in Eq. (\ref{e29}) the LS solution can be obtained using the procedure described in Eqs. (\ref{e06}-\ref{e07}). Equation (\ref{e29}) has been tested, providing excellent results, which are not included for sake of brevity.

\subsection{Constraints: $y (t_1) = y_1$ and $\ddot{y} (t_2) = \ddot{y}_2 \; \to \; \ddot{y}_{2x} = \ddot{y}_2 \, \dfrac{4}{\delta t^2}$}

For this case the constrained equation
\begin{equation*}
    y (x) = g (x) - x \, (y_1 - g_1) + \dfrac{x^2 + x}{2} \, (\ddot{y}_{2x} - \ddot{g}_2)
\end{equation*}
can be used. Substituting this equation in Eq. (\ref{e10}), we obtain
\begin{equation}\label{e30}
\begin{array}{c}
     \dfrac{4}{\delta t^2} f_2 \left(\dfrac{\dd^2 g}{\dd x^2} - \ddot{g}_2\right) +  \dfrac{2}{\delta t} f_1 \left(\dfrac{\dd g}{\dd x} + g_1 - \ddot{g}_2 \dfrac{2 x + 1}{2}\right) + f_0 \left(g + g_1 x - \ddot{g}_2 \, \dfrac{x^2 + x}{2}\right) = \\ = f - \dfrac{4}{\delta t^2} \ddot{y}_{2x} f_2 -  \dfrac{2}{\delta t} f_1 \left(- y_1 + \ddot{y}_{2x} \dfrac{2x + 1}{2}\right) - f_0 \left(- y_1 x + \ddot{y}_{2x} \dfrac{x^2 + x}{2}\right)
\end{array}
\end{equation}
Then, using the expressions provided in Eqs. (\ref{e13}-\ref{e16}) in Eq. (\ref{e30}) the LS solution can be obtained using the procedure described in Eqs. (\ref{e06}-\ref{e07}). Equation (\ref{e30}) has been tested, providing excellent results, which are not included for sake of brevity.

\subsection{Constraints: $\dot{y} (t_1) = \dot{y}_1 \; \to \; \dot{y}_{1x} = \dot{y}_1 \, \dfrac{2}{\delta t}$ and $y (t_2) = y_2$}

For this case the constrained equation
\begin{equation*}
    y (x) = g (x) + (y_2 - g_2)  + (x - 1) (\dot{y}_{1x} - \dot{g}_1)
\end{equation*}
can be used. Substituting this equation in Eq. (\ref{e10}), we obtain
\begin{equation}\label{e31}
    \begin{array}{c} \dfrac{4}{\delta t^2} f_2 \dfrac{\dd^2 g}{\dd x^2} + \dfrac{2}{\delta t} f_1 \left(\dfrac{\dd g}{\dd x} - \dot{g}_1\right) + f_0 \left[g - g_2 - \dot{g}_1 (x - 1)\right] = \\ = f - \dfrac{2}{\delta t} \dot{y}_{1x} f_1 - f_0 \left[y_2 + \dot{y}_{1x} (x - 1)\right]\end{array}
\end{equation}
Then, using the expressions provided in Eqs. (\ref{e13}-\ref{e16}) in Eq. (\ref{e31}) the LS solution can be obtained using the procedure described in Eqs. (\ref{e06}-\ref{e07}). Equation (\ref{e31}) has been tested, providing excellent results, which are not included for sake of brevity.

\subsection{Constraints: $\dot{y} (t_1) = \dot{y}_1 \; \to \; \dot{y}_{1x} = \dot{y}_1 \, \dfrac{2}{\delta t}$ and $\dot{y} (t_2) = \dot{y}_2 \; \to \; \dot{y}_{2x} = \dot{y}_2 \, \dfrac{2}{\delta t}$}

For this case the constrained equation
\begin{equation*}
    y (x) = g (x) + \dfrac{x}{2} \left(1 - \dfrac{x}{2}\right)(\dot{y}_{1x} - \dot{g}_1) + \dfrac{x}{2} \left(1 + \dfrac{x}{2}\right) (\dot{y}_{2x} - \dot{g}_2)
\end{equation*}
can be used. Substituting this equation in Eq. (\ref{e10}), we obtain
\begin{multline}\label{e32}
    \dfrac{4}{\delta t^2} f_2 \left(\dfrac{\dd^2 g}{\dd x^2} + \dfrac{\dot{g}_1 - \dot{g}_2}{2}\right) + \dfrac{2}{\delta t} f_1 \left[\dfrac{\dd g}{\dd x} - \dfrac{\dot{g}_1 (1 - x) + \dot{g}_2 (x + 1)}{2}\right] + \\ + f_0 \left\{g - \dfrac{x}{2} \left[\dot{g}_1 \left(1 - \dfrac{x}{2}\right) + \dot{g}_2 \left(\dfrac{x}{2} + 1\right)\right]\right\} = f - \dfrac{2}{\delta t^2} (\dot{y}_{2x} - \dot{y}_{1x}) f_2 \\ - \dfrac{1}{\delta t} f_1 \left[\dot{y}_{1x} (1 - x) + \dot{y}_{2x} (x + 1)\right] - f_0 \dfrac{x}{2} \left[\dot{y}_{1x} \left(1 - \dfrac{x}{2}\right) + \dot{y}_{2x} \left(\dfrac{x}{2} + 1\right)\right]
\end{multline}
Then, using the expressions provided in Eqs. (\ref{e13}-\ref{e16}) in Eq. (\ref{e32}) the LS solution can be obtained using the procedure described in Eqs. (\ref{e06}-\ref{e07}).
Equation (\ref{e32}) has been tested for the special case of the Mathieu's DE \cite{Mathieu1}.

\subsection{Constraints: $\dot{y} (t_1) = \dot{y}_1 \; \to \; \dot{y}_{1x} = \dot{y}_1 \, \dfrac{2}{\delta t}$ and $\ddot{y} (t_2) = \ddot{y}_2 \; \to \; \ddot{y}_{2x} = \ddot{y}_2 \, \dfrac{4}{\delta t^2}$}

For this case the constrained equation
\begin{equation*}
    y (x) = g(x) + x \, (\dot{y}_{1x} - \dot{g}_1) + x \, \left(\dfrac{x}{2} + 1\right) (\ddot{y}_{2x} - \ddot{g}_2)
\end{equation*}
can be used. Substituting this equation in Eq. (\ref{e10}), we obtain
\begin{equation}\label{e33}
\begin{array}{c}
    \dfrac{4}{\delta t^2} f_2 \left(\dfrac{\dd^2 g}{\dd x^2} - \ddot{g}_2\right) + \dfrac{2}{\delta t} f_1 \left[\dfrac{\dd g}{\dd x} - \dot{g}_1 - \ddot{g}_2 (x + 1)\right] + f_0 \left[g - \dot{g}_1 x - \ddot{g}_2 \left(\dfrac{x}{2} + 1\right) x\right] = \\ = f - \dfrac{4}{\delta t^2} \ddot{y}_{2x} f_2 - \dfrac{2}{\delta t} f_1 \left[\dot{y}_{1x} + \ddot{y}_{2x} (x + 1)\right] - f_0 \left[\dot{y}_{1x} x + \ddot{y}_{2x} \left(\dfrac{x}{2} + 1\right) x\right]
\end{array}
\end{equation}
Then, using the expressions provided in Eqs. (\ref{e13}-\ref{e16}) in Eq. (\ref{e33}) the LS solution can be obtained using the procedure described in Eqs. (\ref{e06}-\ref{e07}). Equation (\ref{e33}) has been tested, providing excellent results, which are not included for sake of brevity.

\subsection{Constraints: $\ddot{y} (t_1) = \ddot{y}_1 \; \to \; \ddot{y}_{1x} = \ddot{y}_1 \, \dfrac{4}{\delta t^2}$ and $y (x_2) = y_2$}

For this case the constrained equation
\begin{equation*}
    y (x) = g (x) + x \, (y_2 - g_2) + \dfrac{x}{2} \, (x - 1) (\ddot{y}_1 - \ddot{g}_1)
\end{equation*}
can be used. Substituting this equation in Eq. (\ref{e10}), we obtain
\begin{equation}\label{e34}
\begin{array}{c}
    \dfrac{4}{\delta t^2} f_2 \left(\dfrac{\dd^2 g}{\dd x^2} - \ddot{g}_1\right) + \dfrac{2}{\delta t} f_1 \left(\dfrac{\dd g}{\dd x} - g_2 - \ddot{g}_1 \dfrac{2 x - 1}{2}\right) + f_0 \left(g -g_2 x - \ddot{g}_1 \, \dfrac{x^2 - x}{2}\right) = \\ = f - \dfrac{4}{\delta t^2} \ddot{y}_1 f_2 - \dfrac{2}{\delta t} f_1 \left(y_2 + \ddot{y}_1 \dfrac{2x - 1}{2}\right) - f_0 \left(y_2 + \ddot{y}_1 \dfrac{x^2 - x}{2}\right)
\end{array}
\end{equation}
Then, using the expressions provided in Eqs. (\ref{e13}-\ref{e16}) in Eq. (\ref{e34}) the LS solution can be obtained using the procedure described in Eqs. (\ref{e06}-\ref{e07}). Equation (\ref{e34}) has been tested, providing excellent results, which are not included for sake of brevity.

\subsection{Constraints: $\ddot{y} (t_1) = \ddot{y}_1 \; \to \; \ddot{y}_{1x} = \ddot{y}_1 \, \dfrac{4}{\delta t^2}$ and $\dot{y} (t_2) = \dot{y}_2 \; \to \; \dot{y}_{2x} = \dot{y}_2 \, \dfrac{2}{\delta t}$}

For this case the constrained equation
\begin{equation*}
    y (x) = g(x) + x \, (\dot{y}_{2x} - \dot{g}_2) + \dfrac{x}{2} \, (x - 2) (\ddot{y}_1 - \ddot{g}_1)
\end{equation*}
can be used. Substituting this equation in Eq. (\ref{e10}), we obtain
\begin{equation}\label{e35}
\begin{array}{c}
    \dfrac{4}{\delta t^2} f_2 \left(\dfrac{\dd^2 g}{\dd x^2} - \ddot{g}_1\right) + \dfrac{2}{\delta t} f_1 \left[\dfrac{\dd g}{\dd x} - \dot{g}_2 - \ddot{g}_1 (x - 1)\right] + f_0 \left[g - \dot{g}_2 x - \ddot{g}_1 \left(\dfrac{x^2}{2} - x\right)\right] = \\ = f - \dfrac{4}{\delta t^2} \ddot{y}_1 f_2 - \dfrac{2}{\delta t} f_1 \left[\dot{y}_{2x} + \ddot{y}_1 (x - 1)\right] - f_0 \left[\dot{y}_{2x} x + \ddot{y}_1 \left(\dfrac{x^2}{2} - x\right)\right]
\end{array}
\end{equation}
Then, using the expressions provided in Eqs. (\ref{e13}-\ref{e16}) in Eq. (\ref{e35}) the LS solution can be obtained using the procedure described in Eqs. (\ref{e06}-\ref{e07}). Equation (\ref{e35}) has been tested, providing excellent results, which are not included for sake of brevity.

\subsection{Constraints: $\ddot{y} (t_1) = \ddot{y}_1 \; \to \; \ddot{y}_{1x} = \ddot{y}_1 \, \dfrac{4}{\delta t^2}$ and $\ddot{y} (t_2) = \ddot{y}_2 \; \to \; \ddot{y}_{2x} = \ddot{y}_2 \, \dfrac{4}{\delta t^2}$}

For this case the constrained equation
\begin{equation*}
    y (x) = g(x) + \dfrac{x^2}{12} \, (3 - x) (\ddot{y}_{1x} - \ddot{g}_1) + \dfrac{x^2}{12} \, (3 + x) (\ddot{y}_{2x} - \ddot{g}_2)
\end{equation*}
can be used. Substituting this equation in Eq. (\ref{e10}), we obtain
\begin{equation}\label{e36}
\begin{array}{c}
    \dfrac{4}{\delta t^2} f_2 \left[\dfrac{\dd^2 g}{\dd x^2} - \dfrac{\ddot{g}_1 (1 - x) + \ddot{g}_2 (x + 1)}{2}\right] + \\ + \dfrac{2}{\delta t} f_1 \left[\dfrac{\dd g}{\dd x} - \dfrac{\ddot{g}_1 (2 x - x^2) + \ddot{g}_2 (x^2 + 2 x)}{4}\right] + \\ + f_0 \left[g - \dfrac{\ddot{g}_1 (3 x^2 - x^3) + \ddot{g}_2 (x^3 + 3 x^2)}{12}\right] = \\ = f - \dfrac{2}{\delta t^2} f_2 [\ddot{y}_{1x} (1 - x) + \ddot{y}_{21x} (x + 1)] - \dfrac{2}{\delta t} f_1 \dfrac{\ddot{y}_1 (2 x - x^2) + \ddot{y}_{2x} (x^2 + 2 x)}{4} \\ - f_0 \dfrac{\ddot{y}_{1x} (3 x^2 - x^3) + \ddot{y}_{2x} (x^3 + 3 x^2)}{12}
\end{array}
\end{equation}
Then, using the expressions provided in Eqs. (\ref{e13}-\ref{e16}) in Eq. (\ref{e36}) the LS solution can be obtained using the procedure described in Eqs. (\ref{e06}-\ref{e07}). Equation (\ref{e36}) has been tested, providing excellent results, which are not included for sake of brevity.

\subsection{Optimal control example: state known at initial time and costate at final time}

This example consists of the linear DE,
\begin{equation}\label{e37}
    \begin{Bmatrix} \dot{\B{x}} \\ \dot{\B{\lambda}}\end{Bmatrix} = \begin{bmatrix} A_{11} (t) & A_{12} (t) \\ A_{21} (t) & A_{22} (t)\end{bmatrix} \begin{Bmatrix} \B{x} \\ \B{\lambda}\end{Bmatrix} \qquad \text{subject to:} \; \left\{\begin{array}{l} \B{x} (t_0) = \B{x}_0 \\ \B{\lambda} (t_f) = \B{\lambda}_f\end{array}\right.,
\end{equation}
with the constrained expressions,
\begin{equation*}
    \left\{\begin{array}{l} \B{x} (t) = \B{g}_x (t) + (\B{x}_0 - \B{g}_{x0}) \\ \B{\lambda} (t) = \B{g}_{\lambda} (t) + (\B{\lambda}_f - \B{g}_{\lambda f})\end{array}\right..
\end{equation*}
Assuming, $\B{x} = \{x, \dot{x}\}\T$ and $\B{\lambda} = \{\lambda_x, \lambda_{\dot{x}}\}\T$, then
\begin{equation*}
    \B{g}_x (t) = \begin{bmatrix} \B{h}\T (t) \\ \dot{\B{h}}\T (t)\end{bmatrix} \B{\alpha} \qquad \text{and} \qquad \B{g}_{\lambda} (t) = \begin{bmatrix} \B{\beta}\T \\ \B{\gamma}\T\end{bmatrix} \B{h} (t)
\end{equation*}
and the constrained expressions become,
\begin{equation*}
    \B{x} (t) = \B{x}_0 + \begin{bmatrix} \B{h}\T - \B{h}_0\T \\ \dot{\B{h}}\T - \dot{\B{h}}_0\T\end{bmatrix} \B{\alpha} \qquad \text{and} \qquad \B{\lambda} (t) = \B{\lambda}_f + \begin{bmatrix} \B{\beta}\T \\ \B{\gamma}\T\end{bmatrix} (\B{h} - \B{h}_f).
\end{equation*}
Then, the dynamic equation becomes,
\begin{equation*}
    \left\{\begin{array}{l}
    \begin{bmatrix} \dot{\B{h}}\T \\ \ddot{\B{h}}\T\end{bmatrix} \B{\alpha} - A_{11} \begin{bmatrix} \B{h}\T - \B{h}_0\T \\ \dot{\B{h}}\T - \dot{\B{h}}_0\T\end{bmatrix} \B{\alpha} - A_{12} \begin{bmatrix} \B{\beta}\T \\ \B{\gamma}\T\end{bmatrix} (\B{h} - \B{h}_f) = A_{11} \B{x}_0 + A_{12} \B{\lambda}_f \\
    \begin{bmatrix} \B{\beta}\T \\ \B{\gamma}\T\end{bmatrix} \dot{\B{h}} - A_{21} \begin{bmatrix} \B{h}\T - \B{h}_0\T \\ \dot{\B{h}}\T - \dot{\B{h}}_0\T\end{bmatrix} \B{\alpha} - A_{22} \begin{bmatrix} \B{\beta}\T \\ \B{\gamma}\T\end{bmatrix} (\B{h} - \B{h}_f) = A_{21} \B{x}_0 + A_{22} \B{\lambda}_f
    \end{array}\right.,
\end{equation*}
which can be written in matrix form,
\begin{equation}\label{e38}
     \mathcal{M} \, \begin{Bmatrix} \B{\alpha} \\ \B{\beta} \\ \B{\gamma}\end{Bmatrix} = \begin{bmatrix} A_{11} & A_{12} \\ A_{21} & A_{22}\end{bmatrix} \begin{Bmatrix} \B{x}_0 \\ \B{\lambda}_f\end{Bmatrix},
\end{equation}
where
\begin{equation*}
     \mathcal{M} = \begin{bmatrix} \begin{bmatrix} \dot{\B{h}}\T \\ \ddot{\B{h}}\T\end{bmatrix} - A_{11} \begin{bmatrix} \B{h}\T - \B{h}_0\T \\ \dot{\B{h}}\T - \dot{\B{h}}_0\T\end{bmatrix} & - A_{12} \begin{bmatrix} \B{h}\T - \B{h}_f\T \\ \B{0}\T\end{bmatrix} & - A_{12} \begin{bmatrix} \B{0}\T \\ \B{h}\T - \B{h}_f\T\end{bmatrix} \\ - A_{21} \begin{bmatrix} \B{h}\T - \B{h}_0\T \\ \dot{\B{h}}\T - \dot{\B{h}}_0\T\end{bmatrix} & \begin{bmatrix} \dot{\B{h}}\T \\ \B{0}\T\end{bmatrix} - A_{22} \begin{bmatrix} \B{h}\T - \B{h}_f\T \\ \B{0}\T\end{bmatrix} & \begin{bmatrix} \B{0}\T \\ \dot{\B{h}}\T\end{bmatrix} - A_{22} \begin{bmatrix} \B{0}\T \\ \B{h}\T - \B{h}_f\T\end{bmatrix} \end{bmatrix}.
\end{equation*}
Equation (\ref{e38}) is linear in the unknown coefficients (vectors: $\B{\alpha}$, $\B{\beta}$, and $\B{\gamma}$) and, therefore, it can be solved by LS as done in the previous numerical examples.

Equation (\ref{e37}) can be subject to different constraints. The following are two examples that can be solved by LS using the corresponding \emph{constrained expressions},
\begin{eqnarray*}
    \left\{\begin{array}{l} \B{x} (t_0) = \B{x}_0 \\ \dot{\B{x}} (t_0) = \dot{\B{x}}_0\end{array}\right. \; \quad & \to & \left\{\begin{array}{l} \B{x} (t) = \B{g}_x (t) + (\B{x}_0 - \B{g}_{x0}) + (t - t_0) (\dot{\B{x}}_0 - \dot{\B{g}}_{x0}) \\ \B{\lambda} (t) = \B{g}_{\lambda} (t)\end{array}\right. \\
    \left\{\begin{array}{l} \B{x} (t_0) = \B{x}_0 \\ \B{\lambda} (t_f) = \B{x} (t_f)\end{array}\right. & \to & \left\{\begin{array}{l} \B{x} (t) = \B{g}_x (t) + (\B{x}_0 - \B{g}_{x0}) \\ \B{\lambda} (t) = \B{g}_{\lambda} (t) + (\B{g}_{xf} + \B{x}_0 - \B{g}_{x0} - \B{g}_{\lambda f})\end{array}\right.
\end{eqnarray*}

\section*{Conclusions and future work}

This study presents a new approach to provide least-squares solutions of linear nonhomogeneous differential equations of any order with nonconstant coefficients, continuous and non singular in the independent variable integration range. For sake of brevity and without loosing generality, the implementation of the proposed method has been applied to second order differential equations. This least-squares approach can be adopted to solve initial and boundary value problems with constraints given in terms of the function and/or derivatives.

The proposed method is based on searching the solution with a specific expression, called \emph{constrained expression}, which is a function with embedded differential equation constraints. This expression is given in terms of a new unknown function, $g (t)$. The original differential equation is rewritten in terms of $g (t)$, thus obtaining a new differential equation where the constraints are embedded in the differential equation itself. Then, the $g (t)$ function is expressed as a linear combination of basis functions, $\B{h} (t)$. The coefficients of this linear expansion are then computed by least-squares by specializing the new differential equation for a set of $N$ different values of the independent variable. In this study the Chebyshev orthogonal polynomial of the first kind have been selected as basis functions. This choice may not be a good choice because each subsequent derivative degree increases the range by approximately one order of magnitude (and because polynomials are in general a bad choice to describe potential periodic solutions).

Numerical tests have been performed for initial value problems with known solution. A direct comparison has been made with the solution provided by MATLAB function \verb"ODE45", implementing the Runge-Kutta-Fehlberg variable-step integrator. In this test, the least-squares approach shows five orders of magnitude accuracy gain. Numerical tests have been performed for boundary value problems for the four cases of known, unknown, no, and infinite solutions. In particular, the condition number and the residual mean of the least-square approach can be used to discriminate whether a boundary value problem has no solution or infinite solutions.

The proposed method is easy to implement, not iterative and, as opposed to classic numerical approaches, the solution error distribution do not increase along the integration but it is approximately uniformly distributed in the integration range. In addition, the proposed technique is identical to solve initial and boundary value problems. The method can also be used to solve higher order linear differential equations with linear constraints.

This study is not complete as many investigations should be performed before setting this approach as a standard way to integrate linear differential equations. Many research areas are still open, full of question marks. A few of these open research areas are:
\begin{enumerate}
  \item Extension to weighted least-squares;
  \item Nonuniform distribution of points and optimal distributions of points to increase accuracy in specific ranges of interest;
  \item Comparisons with different function bases and identification of an optimal function base (if it exists!);
  \item Analysis using Fourier bases;
  \item Accuracy analysis of number of basis functions versus points distribution;
  \item Extension to nonlinear differential equations;
  \item Extension to partial differential equations.
  \item Extension to nonlinear constraints.
\end{enumerate}
This study does not provide answers to the above questions, but it provides suggestions of important area of research and basic tools to dig.


\begin{thebibliography}{9}
\bibitem{Mortari} Mortari, D. ``The Theory of Connections. Part 1: Connecting Points,'' AAS 17-256 of 2017 AAS/AIAA Space Flight Mechanics Meeting Conference, San Antonio, TX, February 5-9, 2017.
\bibitem{Strang} Strang, G. Differential Equations and Linear Algebra, Wellesley-Cambridge Press, 2015, ISBN 0980232791, 9780980232790.
\bibitem{Lin} Lin, Y., Enszer, J.A., and Stadtherr, M.A. ``Enclosing all Solutions of Two-Point Boundary Value Problems for ODEs,'' \emph{Computers and Chemical Engineering}, 2008, pp. 1714-1725.
\bibitem{Venkataraman1} Venkataraman, P. ``A New Class of Analytical Solutions to Nonlinear Boundary Value Problems,'' DETC2005-84604, 25th Computers and Information in Engineering (CIE) Conference, Long Beach, CA, September, 2005.
\bibitem{Venkataraman2}  Venkataraman, P. ``Explicit Solutions for Linear Boundary Value Problems using B\'ezier Functions,'' DETC2006-99227, 26th Computers and Information in Engineering (CIE) Conference, Philadelphia, PA, September, 2006.
\bibitem{Venkataraman3} Venkataraman, P. and Michopoulos, J.G. ``Explicit Solutions for Nonlinear Partial Differential Equations,'' DETC2007-35439, 27th Computers and Information in Engineering (CIE) Conference, Las Vegas, NA, September, 2007.
\bibitem{Zheng} Zheng, J.M., Sederberg, T.W., and Johnson, R.W. ``Least-Squares Methods for Solving Differential Equations using B\'ezier Control Points,'' \emph{Applied Numerical Mathematics}, Vol. 48, No. 2, Feb. 2004, pp. 237-252.
\bibitem{Evrenosoglu} Evrenosoglu, M., Somali, S. ``Least-Squares Methods for Solving Singular Perturbed Two-Point Boundary Value Problems Using B\'ezier Control Points,'' \emph{Applied Mathematics Letters}, Vol. 21, No. 10, Oct. 2008, pp. 1029–-1032.
\bibitem{Ghomanjani} Ghomanjani, F., Kamyad, A.V., and Kiliçman, A. ``B\'ezier Curves Method for Fourth-Order Integro-Differential Equations,'' \emph{Abstract and Applied Analysis}, Vol. 2013, Article ID 672058, 5 pages, 2013. doi:10.1155/2013/672058.
\bibitem{Doha} Doha, E.H., Bhrawy, A.H., and Saker, M.A. ``On the Derivatives of Bernstein Polynomials: An Application for the Solution of High Even-Order Differential Equations,'' Hindawi Publishing Corporation, Boundary Value Problems, Vol. 2011, Article ID 829543, 16 pages, doi:10.1155/2011/829543
\bibitem{Mathieu1} Mathieu, E. ``M\'emoire sur Le Mouvement Vibratoire d'une Membrane de Forme Elliptique,'' \emph{Journal de Math\'ematiques Pures et Appliqu\'ees}, 137-203, 1868.
\bibitem{Mathieu2} Meixner, J. and Sch\"{a}fke, F.W. ``Mathieusche Funktionen und Sph\"{a}roidfunktionen'' Springer, 1954.
\end{thebibliography}
\end{document}